\documentclass[12pt]{amsart}
\usepackage{amssymb,mathrsfs}
\newtheorem{theorem}{Theorem}[section]
\newtheorem{lemma}[theorem]{Lemma}

\newtheorem{remark}[theorem]{Remark}
\newtheorem{definition}[theorem]{Definition}
\newtheorem{conjecture}[theorem]{Conjecture}
\newtheorem{example}[theorem]{Example}

\newtheorem{problem}[theorem]{Problem}

\renewcommand{\theequation}{\thesection.\arabic{equation}}
\renewcommand{\thetheorem}{\thesection.\arabic{theorem}}

\title{The Modularity of an Abelian Variety}

\begin{document}

\author{Jae-Hyun Yang}

\address{Yang Institute for Advanced Study
\newline\indent
Seoul 07989, Korea
\vskip 2mm
and
\vskip 2mm
Department of Mathematics
\newline\indent
Inha University
\newline\indent
Incheon 22212, Korea}

\email{jhyang@inha.ac.kr\ \ or\ \ yangsiegel@naver.com}

\renewcommand{\theequation}{\thesection.\arabic{equation}}
\renewcommand{\thetheorem}{\thesection.\arabic{theorem}}
\renewcommand{\thelemma}{\thesection.\arabic{lemma}}
\newcommand{\BR}{\mathbb R}
\newcommand{\BQ}{\mathbb Q}
\newcommand{\BF}{\mathbb F}
\newcommand{\BT}{\mathbb T}
\newcommand{\BM}{\mathbb M}
\newcommand{\bn}{\bf n}
\def\charf {\mbox{{\text 1}\kern-.24em {\text l}}}
\newcommand{\BC}{\mathbb C}
\newcommand{\BZ}{\mathbb Z}

\thanks{
\newline \indent
Siegel Hecke eigenforms, zeta functions, Hecke algebra, $p$-Satake parameters.}

\begin{abstract}
We introduce the concept of the modularity of an abelian variety defined over the rational number field
extending the modularity of an elliptic curve. We discuss the modularity of an abelian variety over $\BQ$. We conjecture that a simple abelian variety over $\BQ$ is modular.
\end{abstract}

\maketitle

\newcommand\tr{\triangleright}
\newcommand\al{\alpha}
\newcommand\be{\beta}
\newcommand\g{\gamma}
\newcommand\gh{\Cal G^J}
\newcommand\G{\Gamma}
\newcommand\de{\delta}
\newcommand\e{\epsilon}
\newcommand\z{\zeta}
\newcommand\vth{\vartheta}
\newcommand\vp{\varphi}
\newcommand\om{\omega}
\newcommand\p{\pi}
\newcommand\la{\lambda}
\newcommand\lb{\lbrace}
\newcommand\lk{\lbrack}
\newcommand\rb{\rbrace}
\newcommand\rk{\rbrack}
\newcommand\s{\sigma}
\newcommand\w{\wedge}
\newcommand\fgj{{\frak g}^J}
\newcommand\lrt{\longrightarrow}
\newcommand\lmt{\longmapsto}
\newcommand\lmk{(\lambda,\mu,\kappa)}
\newcommand\Om{\Omega}
\newcommand\ka{\kappa}
\newcommand\ba{\backslash}
\newcommand\ph{\phi}
\newcommand\M{{\Cal M}}
\newcommand\bA{\bold A}
\newcommand\bH{\bold H}
\newcommand\D{\Delta}

\newcommand\Hom{\text{Hom}}
\newcommand\cP{\Cal P}

\newcommand\cH{\Cal H}

\newcommand\pa{\partial}

\newcommand\pis{\pi i \sigma}
\newcommand\sd{\,\,{\vartriangleright}\kern -1.0ex{<}\,}
\newcommand\wt{\widetilde}
\newcommand\fg{\frak g}
\newcommand\fk{\frak k}
\newcommand\fp{\frak p}
\newcommand\fs{\frak s}
\newcommand\fh{\frak h}
\newcommand\Cal{\mathcal}

\newcommand\fn{{\frak n}}
\newcommand\fa{{\frak a}}
\newcommand\fm{{\frak m}}
\newcommand\fq{{\frak q}}
\newcommand\CP{{\mathcal P}_n}
\newcommand\Hnm{{\mathbb H}_n \times {\mathbb C}^{(m,n)}}
\newcommand\BD{\mathbb D}
\newcommand\BH{\mathbb H}
\newcommand\CCF{{\mathcal F}_n}
\newcommand\CM{{\mathcal M}}
\newcommand\Gnm{\Gamma_{n,m}}
\newcommand\Cmn{{\mathbb C}^{(m,n)}}
\newcommand\Yd{{{\partial}\over {\partial Y}}}
\newcommand\Vd{{{\partial}\over {\partial V}}}

\newcommand\Ys{Y^{\ast}}
\newcommand\Vs{V^{\ast}}
\newcommand\LO{L_{\Omega}}
\newcommand\fac{{\frak a}_{\mathbb C}^{\ast}}

\vskip 3mm

\begin{section}{{\bf Introduction}}
\setcounter{equation}{0}
\vskip 3mm
An elliptic curve $E$ over $\BQ$ is said to be ${\sf modular}$
if the $L$-function $L(E,s)$ of $E$ equals the $L$-function $L(f,s)$ for
some eigenform $f$, equivalently if $E$ has a finite covering by a modular curve
of the form $X_0(N)$. At the Tokyo-Nikko conference held in 1955, Yutaka Taniyama\,(1927--1958)
made a suggestion that every elliptic curve over $\BQ$ is modular. At that time
his suggestion was not clear and hence was not known in the mathematics community.
In the early 1960's, Goro Shimura\,(1930--2019) refined Taniyama's suggestion through private
conversations with a number of mathematicians and their efforts. In particular
he discussed this subject with Andr{\'e} Weil\,(1906--1998) seriously and intensively.
Weil gave conceptual evidence for Taniyama's suggestion in his famous paper
\cite{Weil} published in 1967. Through Weil's paper, this suggestion was widely known
as the so-called {\sf Shimura-Taniyama conjecture} in the mathematics community.

\vskip 2mm
The Shimura-Taniyama conjecture associates objects of representation theory to
those of algebraic geometry. It states that the $L$-series of an elliptic curve over $\BQ$
which measures the behaviour of the curve mod $p$ for all primes $p$, can be identified
with an integral transform of the Fourier series defined from an eigenform.

\vskip 2mm
In 1985 Frey \cite{Frey} made the remarkable observation that the Shimura-Taniyama conjecture
for a semistable elliptic curve over $\BQ$ should imply Fermat's Last Theorem.
The precise mechanism relating the two was formulated by Serre \cite{S} as
the $\epsilon$-conjecture and this conjecture was proved by Ribet \cite{Ribet} in the summer of 1986.
Ribet's result only requires to prove the Shimura-Taniyama conjecture for a semistable
elliptic curve over $\BQ$ in order to reduce Fermat's Last Theorem. As soon as Andrew Wiles
learned Ribet's result, he began to work the Shimura-Taniyama conjecture for a semistable
elliptic curve over $\BQ$ in the late summer of 1986.
Here the semistabilty of an elliptic curve $E$ over $\BQ$ is defined as follows\,:

\begin{definition}\label{def:1.1}
An elliptic curve $E$ over $\BQ$ is said to be {\sf semistable} at the prime $q$
if it is isomorphic to an elliptic curve $\tilde E$ over $\BQ$ which
$\tilde E$\,{\rm (mod} $q${\rm )} is either nonsingular or has a node. An elliptic curve
over $\BQ$ is called {\sf semistable} if it is semistable at every prime.
\end{definition}

\vskip 2mm
On June 21-23 in 1993, Wiles had given a series of lectures under the title,
''{\sf Elliptic curves, modular forms and Galois representations}" at the
Issac Newton Institute for Mathematical Sciences in Cambridge, England.
In this last lecture on June 23, Wiles commented that he had proved a part of
the Shimura-Taniyama conjecture to the effect that every semistable elliptic
curve over $\BQ$ is modular. With the aide of the works of Frey \cite{Frey}, Serre \cite{S}
and Ribet \cite{Ribet}, his proof solves Fermat's Last Theorem which had been unsolved
for more than 350 years. The news spread out all over the world through the well-known
newspapers and magazines because Fermat's Last Theorem holds great fascination for
amateurs and professionals alike. But in the fall that year it turned out that
the proof of Wiles was incomplete and flawed. Precisely, his construction of the Euler
system used to extend Flach's method was not complete. He was very concerned with
filling the gaps in his flawed proof. At that time he was completely isolated from outside.
In January 1994, he proposed Dr. Richard Taylor, his former Ph.D. student in
Cambridge University, UK to join him
in the attempt to repair the Euler system argument. Wiles still was convinced that his
method was correct. Dr. Taylor accepted his proposal and joined in that project.
On September 19th in 1994, Wiles was quite convinced that his method was correct and his gap
could be filled up. After he invited Dr. Taylor to Princeton again, he completed his proof
with the aid of Dr. Taylor on October, 1994 and submitted his paper to Annals of Mathematics
on October 14, 1994. Finally his paper was accepted and was published in May, 1995
(cf.\,\cite{W, TW}). We refer to \cite{Mozzo} for more interesting stories.

\vskip 2mm
The aim of this article is to introduce the concept of the modularity of an abelian variety
defined over the rational number field extending the modularity of an elliptic curve and to
discuss the modularity of an abelian variety over $\BQ$. An abelian variety is defined as
follows in \cite[p.\,39]{MU}\,:

\begin{definition}\label{def:1.2}
An abelian variety $X$ is a complete\footnote{This means, in particular, that it is irreducible.
See the footnote in \cite[p.\,39]{MU}.}
algebraic variety over an algebraically closed field
$K$ with a group law $m_X:X\times X \lrt X$ such that both $m_X$ and the inverse map are
morphisms of varieties. An abelian variety over $\BQ$  is a connected and complex projective
manifold that is also a group variety generalizing elliptic curves over $\BQ$.
\end{definition}

\vskip 2mm
The paper is organized as follows. In Section 2, we briefly outline the modularity of an
elliptic curve over $\BQ$. Wiles proved the modularity of a semistable elliptic curve over
$\BQ$ leading to solve Fermat's Last Theorem (cf.\,\cite{W, TW}). In 2001, Breuil, Conrad, Diamond
and Taylor proved that every elliptic curve over $\BQ$ is modular (cf.\,\cite{BCDT}).
In Section 3, we introduce the notion of the modularity of an abelian variety over $\BQ$.
We propose several conjectures and open problems.
In Section 4 (Appendix), we describe and survey the Hecke algebra for a symplectic group and
$p$-Satake parameters for the readers.

\vskip 0.51cm \noindent {\bf Notations:} \ \ We denote by
$\BQ,\,\BR$ and $\BC$ the field of rational numbers, the field of
real numbers and the field of complex numbers respectively. We
denote by $\BZ$ and $\BZ^+$ the ring of integers and the set of
all positive integers respectively. For a prime $p$, we denote by $\BQ_p$
the field of $p$-adic numbers and by $\BZ_p$ the ring of $p$-adic integers.
$\overline{\BQ}$\,(resp. $\overline{\BQ_p}$)  denotes the
algebraic closure of $\BQ$\,(resp. $\BQ_p$).
$G_\BQ={\rm Gal}(\overline \BQ/\BQ)$ denotes the absolute Galois group of $\BQ$.
The symbol ``:='' means that
the expression on the right is the definition of that on the left.
For two positive integers $k$ and $l$, $F^{(k,l)}$ denotes the set
of all $k\times l$ matrices with entries in a commutative ring
$F$.
For a square matrix $M\in F^{(k,k)}$ of degree $k$,
${{\rm Tr}(M)}$ denotes the trace of $M$.
For any $M\in F^{(k,l)},\
^t\!M$ denotes the transpose of a matrix $M$. $I_n$ denotes the
identity matrix of degree $n$. For $A\in F^{(k,l)}$ and $B\in
F^{(k,k)}$, we set $B[A]=\,^tABA$ (Siegel's notation).
For a number field $F$, we denote
by ${\mathbb A}_F$ the ring of adeles of $F$. If $F=\BQ$, the
subscript will be omitted.
\begin{equation*}
  Sp(2g,\BR)=\left\{ \,M\in Sp(2g,\BR)\,|\ {}^tMJ_gM=J_g[M]=J_g\,\right\}
\end{equation*}
denotes the symplectic group of degree $g$, where
$$
J_g=\begin{pmatrix}
        0 & I_g \\
        -I_g & 0
      \end{pmatrix}.
$$
\begin{equation*}
  \Gamma_g:=Sp(2g,\BZ)=\left\{ \
  \begin{pmatrix}
        A & B \\
        C & D\end{pmatrix}\in Sp(2g,\BR)\Big| \ A,B,C,D\ \rm{integral}\,\right\}
\end{equation*}
denotes the Siegel modular group of degree $g$.
For a prime $\ell$ and a positive integer $n$, ${\mathbb F}_{\ell^n}$ denotes the field with $\ell^n$ elements. The algebraic closure $\overline{\mathbb F}_\ell$ of ${\mathbb F}_\ell$ is given by
$$\overline{\mathbb F}_\ell=\bigcup_{n\geq 1}{\mathbb F}_{\ell^n}.$$
If $p$ is prime and $\lambda$ is a prime ideal dividing $p$ in the ring of integers in
$\overline \BQ$, there exist a filtration
$$ I_\lambda\subset D_\lambda\subset G_\BQ,$$
where the decomposition group $D_\lambda$ and the inertia group $I_\lambda$ are defined
respectively by
$$ D_\lambda:=\{\,\sigma\in G_\BQ\,|\,\sigma(\lambda)=\lambda\,\} $$
and
$$
I_\lambda:=\{\,\sigma\in D_\lambda\,|\,\sigma(x)\equiv x\,({\rm mod}\,\lambda)\
{\rm{for\ all\ algebraic\ integers}}\ x\,\}.
$$
Then there are natural identifications
$$
D_\lambda\cong {\rm Gal} (\overline{\BQ}_p/\BQ_p)\quad {\rm and}\quad
D_\lambda/I_\lambda\cong {\rm Gal} (\overline{\BF}_p/\BF_p).
$$
We denote by ${\rm Frob}_\lambda\in D_\lambda/I_\lambda$ the inverse image of the
canonical generator $x\mapsto x^p$ of ${\rm Gal}(\overline{\mathbb F}_p/{\mathbb F}_p).$
If $\lambda'$ is another prime lying above $p$, then $\lambda'=\sigma(\lambda)$ for
some $\sigma\in G_\BQ$ and
$$
D_{\lambda'}=\sigma D_\lambda \sigma^{-1},\quad I_{\lambda'}=\sigma I_\lambda \sigma^{-1}
\quad {\rm and}\quad {\rm Frob}_{\lambda'}=\sigma\, {\rm Frob}_{\lambda}\, \sigma^{-1}.
$$
Since we will care about these objects only up to conjugation, we will write
$D_\lambda$ and $I_\lambda$. Now we will write ${\rm Frob}_p$ for any representative of
a ${\rm Frob}_{\lambda}$. If $\rho$ is a representation of $G_\BQ$ which is unramified
at $p$, then ${\rm Tr}(\rho({\rm Frob}_p))$ and $\det (\rho({\rm Frob}_p))$ are well
defined, that is, are independent of the choice of $\lambda$.

\end{section}

\vskip 10mm

\vskip 3mm

\begin{section}{{\bf The Modularity of an Elliptic Curve}}
\setcounter{equation}{0}
\vskip 3mm
We set $\G_1:=SL(2,\BZ)$. For a positive integer $N$, we let $\G (N),\ \G_1 (N)$ and
$\G_0 (N)$ be the congruence subgroups of $\G_1$ such that $\G (N)\subset \G_1(N)
\subset \G_0 (N)\subset \G_1.$ We refer to \cite[pp.\,13--14,\,p.\,21]{Di-S}
for the precise definitions and properties of $\G (N),\ \G_1 (N)$ and $\G_0 (N)$.
Let $\BH_1$ be the Poincar{\'e} upper half plance. The quotient
\begin{equation*}
  Y_1(N):=\G_1(N)\ba \BH_1 \ ({\rm resp}.\  Y_0(N):=\G_0(N)\ba \BH_1)
\end{equation*}
be the complex manifold which has a natural model $Y_1(N)/\BQ$ (resp. $Y_0(N)/\BQ$).
We let $X_1(N)$ (resp. $X_0(N)$) be the smooth projective curve which contains
$Y_1(N)$ (resp. $Y_0(N)$) as a dense Zariski open subset (cf.\,see \cite[pp.\,45--60]{Di-S}).

\vskip 3mm
Let $S_k(N)$ be the space of cusp forms of weight $k\geq 1$ and level $N\geq 1$.
Here $k$ and $N$ be positive integers. We recall that
if $f\in S_k(N)$, it satisfies the following properties\,:
\vskip 2mm
(C1) $f((a\tau+b)(c\tau+d)^{-1})=(c\tau+d)^k f(\tau)$ for all
$\begin{pmatrix}
   a & b \\
   c & d
 \end{pmatrix}\in \G_1(N)$ and $\tau\in \BH_1$;
\vskip 2mm
(C2) $| f(\tau)|^2 ({\rm Im}\,\tau)^k$ is bounded in $\BH_1$ and
\vskip 2mm
(C3) the Fourier expansion of $f(\tau)$ is given by
$$ f(\tau)=\sum_{n=1}^{\infty} a_n (f)\, q^n,\quad {\rm where}\ q=e^{2\pi i\tau}.$$

We define the $L$-series of $f\in S_k(N)$ to be
$$ L(f,s):=\sum_{n=1}^{\infty} a_n (f)\,n^{-s}.$$

\vskip 2mm
For each prime $p\!\not| N$, we recall that the Hecke operator $T_p:S_k(N)\lrt S_k(N)$ is
defined by
\begin{equation*}
  (T_pf)(\tau)=p^{-1}\sum_{i=0}^{p-1} f\left(\frac{\tau+i}{p}\right)
  +p^{k-1} f\left(\frac{ap\tau+b}{cp\tau+d}\right), \quad f\in S_k (N)
\end{equation*}
for any $\begin{pmatrix}
   a & b \\
   c & d
 \end{pmatrix}\in \G_1$ with $c\equiv 0\,({\rm mod}\,N)$ and $d\equiv p\,({\rm mod}\,N)$.
We refer to \cite[p.\,844]{BCDT} or \cite[pp.\,170--171]{Di-S} for more details. The Hecke
operators $T_p\,(p\!\!\not| N)$ can be simultaneously diagonalized on $S_k (N)$ and
a simultaneous eigenvector a {\it Hecke\ eigenform} or simply an {\sf eigenform}.

\vskip 2mm
Let $\lambda$ be a place of the algebraic closure $\bar \BQ$ of $\BQ$ in $\BC$ lying over
a rational integer $\ell$ and ${\bar \BQ}_\lambda$ denote the algebraic closure of $\BQ_\ell$
via $\lambda$. Let $G_\BQ:={\rm Gal}({\bar \BQ}/\BQ)$ be the absolute Galois group of $\BQ$.
It is well known that if $f\in S_k (N)$ is a normalized eigenform with
$a_1(f)=1$, then there exists a unique continuous irreducible Galois representation
\begin{equation*}
  \rho_{f,\lambda}:G_\BQ\lrt GL(2,{\bar \BQ}_\lambda)
\end{equation*}
such that $\rho_{f,\lambda}$ is unramified at $p$ for all primes $p\!\!\not| \,\ell N$ and
\begin{equation*}
  {\rm Tr}\left(\rho_{f,\lambda}({\rm Frob}_p)\right)= a_p(f)\quad {\rm for\ any\ prime}\
p\!\!\not| \, \ell N.
\end{equation*}
The existence of $\rho_{f,\lambda}$ is due to Shimura if $k=2$ \cite{Sh2}, due to Deligne if
$k>2$ \cite{De} and due to Deligne and Serre if $k=1$ \cite{D-S}. We see that $\rho_{f,\lambda}$ is
odd in the sense that $\det \rho_{f,\lambda}$ of complex conjugation is -1. Moreover
$\rho_{f,\lambda}$ is potentially semi-stable at $\ell$ in the sense of Fontaine \cite{F-M}.

\vskip 3mm
We may choose a conjugate of $\rho_{f,\lambda}$ which is valued in
$GL(2,{\mathcal O}_{\bar{\BQ}_\lambda})$ and reducing modulo the maximal ideal and
semi-simplyfing yields an irreducible continuous representation
\begin{equation*}
{\overline \rho}_{f,\lambda}:G_\BQ\lrt GL(2,\overline{\mathbb F}_\ell)
\end{equation*}
which, up to isomorphism, does not depend on the choice of conjugate of $\rho_{f,\lambda}$.

\begin{definition}\label{def:1.1}
Let
$$\rho: G_\BQ\lrt GL(2,\bar{\BQ}_\ell)$$
be an irreducible continuous Galois representation which is unramified outside finitely
many primes and for which the restriction of $\rho$ to a decomposition group at $\ell$ is
potentially semi-stable at $\ell$ in the sense of Fontaine. Then $\rho$ is called
${\sf modular}$ if $\rho$ is equivalent to $\rho_{f,\lambda}$
(denoted $\rho\sim \rho_{f,\lambda}$) for some normalized eigenform $f$ and some place
$\lambda | \ell.$
\end{definition}

\begin{definition}\label{1.2}
Let
$${\overline\rho}: G_\BQ\lrt GL(2,\overline{\mathbb F}_\ell)$$
be a two-dimensional irreducible continuous representation of $G_\BQ$.
Then ${\overline\rho}$ is called
${\sf modular}$ if ${\overline\rho}\sim {\overline\rho}_{f,\lambda}$
for some normalized eigenform $f$ and some place $\lambda | \ell.$
\end{definition}

\vskip 3mm
Let $E$ be an elliptic curve over $\BQ$. We define
\begin{equation*}
  a_p(E):=p+1-|E({\mathbb F}_p)| \quad {\rm for\ a\ prime}\ p.
\end{equation*}
The $L$-function $L(E,s)$ of $E$ is defined by the product of the local $L$-factors
\begin{equation*}
  L(E,s):=\prod_{p|D}\left( \frac{1}{1-a_p(E)p^{-s}}\right)
  \prod_{p\not| D}\left( \frac{1}{1-a_p(E)\,p^{-s}+p^{1-2s}}\right).
\end{equation*}
Then $L(E,s)$ converges absolutely for ${\rm Re}\,s> \frac{3}{2}$ and extends to
an entire function by \cite{BCDT}.
\begin{definition}\label{1.3}
An elliptic curve $E$ over $\BQ$ is called
${\sf modular}$ if there exists a Hecke eigenform $f\in S_2(N)$ such that
\begin{equation*}
  L(E,s)=L(f,s).
\end{equation*}
\end{definition}
Let $E$ be an elliptic curve over $\BQ$ with its conductor $N(E)$.
Let $$\rho_{E,\ell}:G_\BQ \lrt GL(2,\bar{\BQ}_\ell)$$
be the $\ell$-adic representation of $G_\BQ$ with the Tate module
$T_\ell (E)$ as its representation space. Let
\begin{equation*}
  J_1(N):=\Omega^1(X_1(N))^{\vee}/H_1(X_1(N),\BZ)\cong S_2(\G_1(N))^{\vee}/H_1(X_1(N),\BZ)
\end{equation*}
be the Jacobian variety of the modular curve $X_1(N)$. Here $\Omega^1(X_1(N))$
denotes the complex vector space of holomorphic 1-forms on $X_1(N)$ and $W^{\vee}$
denotes the dual space of a complex vector space $W$.
It is known that the following statements are equivalent\,:
\vskip 2mm
(a) $E$ is modular.
\vskip 2mm
(b) There is a non-constant holomorphic mapping $X_1(N)\lrt E(\BC)$ for some
pos- \par
\ \ \ \ \ \,itive integer $N$.
\vskip 2mm
(c) There is a non-constant holomorphic mapping $J_1(N)\lrt E(\BC)$ for some pos- \par
\ \ \ \ \ \,itive integer $N$.
\vskip 2mm
(d) $\rho_{E,\ell}$ is modular for a prime $\ell$.

\vskip 3mm
The above statements have been called the Taniyama-Shimura conjecture. We refer to
\cite{Sh3} for the historical story of this conjecture.
The implication (a)\,$\Longrightarrow$\,(b) follows from a construction of Shimura\,\cite{Sh2}
and a theorem of Faltings\,\cite{F}. The implication (b)\,$\Longrightarrow$\,(d) is due to
Mazur\,\cite{MA}. The implication (d)\,$\Longrightarrow$\,(a) follows from a theorem of
Carayol\,\cite{C}. The implication (c)\,$\Longrightarrow$\,(b) is obvious. Wiles\,\cite{W,TW}
proved that a semistable elliptic curve over $\BQ$ is modular by proving the statement (d).
Thereafter Breuil, Conrad, Diamond and Taylor\,\cite{BCDT} proved that every elliptic curve
over $\BQ$ is modular.

\vskip 3mm
Serre \cite{S} conjectured the following\,:
\vskip 2mm\noindent
{\bf Serre's\ Modularity\ Conjecture\,:} Let $\overline{\rho}:G_\BQ\lrt GL(2,{\mathbb F})$ be a two-dimensional
absolutely irreducible, continuous, odd representation of $G_\BQ$. Here $\mathbb F$ is a finite
field of characteristic $p$. Then $\overline{\rho}$ is modular, i.e., arises from (with respect to
some fixed embedding $\imath:\overline{\BQ}\hookrightarrow \overline{\BQ_p})$ a newform $f$ of
some weight $k\geq 2$ and level $N$ prime to $p$.

\vskip 3mm
In 2009, Khare and Wintenberger \cite{K-W1, K-W2} proved that Serre's Modularity Conjecture is true.

\end{section}

\vskip 10mm

\begin{section}{{\bf The Modularity of an Abelian Variety}}
\setcounter{equation}{0}
\vskip 2mm
Let $G:=Sp(2g,\BR)$ and $K=U(g).$
Let
$$\BH_g:=\{\,\Omega\in\BC^{(g,g)}\,|\ \Omega=\,{}^t\Omega,\ \,{\rm{Im}}\,\Omega>0\ \}$$
be the Siegel upper half plane of degree $g$.
Then $G$ acts on ${\mathbb H}_g$ transitively by

\begin{equation}
  \alpha\cdot \Omega=(A\Omega+B)(C\Omega+D)^{-1},
\end{equation}
where $\alpha=\begin{pmatrix}
           A & B \\
           C & D
         \end{pmatrix}\in G$ and $\Omega\in \BH_g.$
The stabilizer of the action (3.1) at $iI_g$ is

\begin{equation*}
  \left\{ \begin{pmatrix} \,A & B \\ -B & A \end{pmatrix} \Big| \ A+iB\in U(g)\,\right\}
  \cong U(g).
\end{equation*}
Thus we get the biholomorphic map
\begin{equation*}
G/K \lrt \BH_g, \qquad \alpha K \longmapsto \alpha\!\cdot\! iI_g,  \quad \alpha\in G.
\end{equation*}
It is known that $\BH_g$ is an Einstein-K{\"a}hler Hermitian symmetric space.

\vskip 3mm
Let $\G_g:=Sp(2g,\BZ)$ be the Siegel modular group of degree $g$.
For a positive integer $N$, we let
\begin{equation*}
  \G_g(N):=\left\{ \g\in \G_g\,|\ \g\equiv I_{2g}\ ({\rm{mod}}\,N)\,\right\}
\end{equation*}
be the the principal congruence subgroup of $\G_g$ of level $N$. Let
\begin{equation*}
  \G_{g,0}(N):=\left\{ \g\in \G_g\,\big|\ \g=\begin{pmatrix}
                                       A & B \\
                                       C & D
                                     \end{pmatrix},\quad C\equiv 0\ ({\rm{mod}}\,N)\,\right\}
\end{equation*}
and
\begin{equation*}
  \G_{g,1}(N):=\left\{ \g\in \G_g\,\big| \
  \g=\begin{pmatrix}
      A & B \\
      C & D
     \end{pmatrix}\equiv \begin{pmatrix}
      I_g & * \\
      0 & I_g
     \end{pmatrix}\ ({\rm{mod}}\,N)\,\right\}
\end{equation*}
be the congruence subgroups of Level $N$. Then we have the relation
\begin{equation*}
  \G_g (N)\subset \G_{g,1} (N)\subset \G_{g,0} (N)\subset \G_g.
\end{equation*}

\begin{definition}\label{def (3.1)}
Let $\G$ be a congruence subgroup of $\G_g$ and let $k$ be a nonnegative integer $k$.
A function $F:\BH_g\lrt \BC$ is called a
$\sf{Siegel\ modular\ form}$ of degree $g$ and weight $k$ with respect to $\G$
if it satisfies the following conditions\,:
\vskip 2mm
{\sf{(S1)}} $F(\Om)$ is holomorphic on $\BH_g$\,;
\vskip 2mm
{\sf{(S2)}} $F(\g\cdot \Om)=(C\Om+D)^k F(\Om)$ \ \ for all
$\g=\begin{pmatrix}
      A & B \\
      C & D
     \end{pmatrix}\in \G$ and $\Om\in \BH_g$\,;
\vskip 2mm
{\sf{(S3)}} $F(\Om)$ is bounded in any domain $Y\geq Y_0 > 0$ in the case $g=1$.
\end{definition}
We denote the space of all Siegel modular forms of degree $g$ and weight $k$
with respect to $\G$ by $[\G,k]$.
\vskip 2mm
We define the so-called $\sf{Siegel\ operator}$
\begin{equation*}
  \Phi_g: [\G_g,k]\lrt [\G_{g-1},k]
\end{equation*}
by
\begin{equation*}
  (\Phi_g(F))(\Om_1):=\lim_{t\lrt\infty}
  F  \begin{pmatrix}
      \Om_1 & 0 \\
      0 & i\,t
     \end{pmatrix},\quad \Om_1\in \BH_{g-1}.
\end{equation*}
Then $\Phi_g$ is a well-defined linear mapping (cf,\,\cite[pp.\,187--189]{M}).
A Siegel modular form $F\in [\G_g,k]$ is called a $\sf{Siegel\ cusp\ form}$ if
$\Phi_g(F)=0$\,(cf.\,\cite[p.\,198]{M}). We denote by $[\G_g,k]_0$
the space of all Siegel cusp forms in $[\G_g,k]$.

\vskip 3mm
Let $\G$ be a congruence subgroup of $\G_g$. If $F\in [\G,k]$, then $F$ has
a Fourier expansion
\begin{equation*}
  F(\Om)=\sum_{T} a(T;F)\,e^{2\pi i\,{\rm Tr}(T\Om)},
\end{equation*}
where $T$ runs through all $g\times g$ half-integral semi-positive symmetric matrices.
Here ${\rm Tr}(M)$ denotes the trace of a square matrix $M$.
Following Hecke's method, Maass\,\cite[pp.\,202--210]{M} associated with
$F(\Om)$ the Dirichlet series
\begin{equation*}
D(F,s):=\sum_{\{ T \}}\frac{a(T;F)}{\varepsilon (T)}\,(\det T)^{-s},
\end{equation*}
where the summation indicates that $T$ runs through a complete set of representatives
of the sets
$$
\left\{ \, T[U]\,|\ U \ \rm{unimodular}\,\right\}, \ T>0
$$
and $\varepsilon (T)$ denotes the number of unimodular matrices $U$ which satisfy
the equation $T[U]=T.$ We note that the numbers $\varepsilon (T)$ are finite.

\begin{definition}\label{def:3.2}
Let $F$ be a nonzero Siegel Hecke eigenform in $[\G_g,k]_0$. Let
$\alpha_{p,0},\alpha_{p,1},\cdots,$\\
\noindent
$\alpha_{p,g}$ be the $p$-Satake parameters of $F$ at
a prime $p$ (cf.\,see Section 4). We define the $\sf{local\ spinor\ zeta\ function}$
$Z_{F,p}(t)$ of $F$ at $p$ by
$$
Z_{F,p}(t):=(1-\alpha_{p,0}\,t)\sum_{r=1}^g\sum_{1\leq i_1<\cdots < i_r\leq g}
(1-\alpha_{p,0}\alpha_{p,i_1}\cdots\alpha_{p,i_r}\,t).
$$
The $\sf{spinor\ zeta\ function}$ $Z_{F}(s)$ of $F$ is defined to be the following function
\begin{equation*}
Z_{F}(s):=\sum_{p:{\rm prime}} Z_{F,p}(p^{-s})^{-1},\quad {\rm Re}\,s \gg 0.
\end{equation*}
Secondly one has the so-called $\sf{standard\ zeta\ function}$ $D_{F}(s)$ of a Siegel Hecke
eigenform $F$ in $[\G_g,k]_0$ defined by
\begin{equation*}
D_{F}(s):=\sum_{p:{\rm prime}} D_{F,p}(p^{-s})^{-1},\quad {\rm Re}\,s \gg 0,
\end{equation*}
where
$$
D_{F,p}(t)=(1-t) \sum_{i=1}^g (1-\alpha_{p,i}\,t)(1-\alpha_{p,i}^{-1}\,t).
$$
We refer to \cite[p.\,249]{Y}.
\end{definition}

\begin{remark}\label{rk:3.3}
(1) If $g=1$, the spinor zeta function $Z_f(s)$ of a Hecke eigenform $f$ in $S_k(\G_1)$
is nothing but the Hecke $L$-function $L(f,s)$ of $f$.
\vskip 2mm\noindent
(2) If $g=1$, the standard zeta function $D_f(s)$ of a Hecke eigenform
$f(\tau)=\sum_{n=1}^{\infty} a(n)\,e^{2\pi i n\tau}$ in $S_k(\G_1)$ has the following equation
$$
D_f(s-k+1)=\sum_{p:{\rm prime}}(1+p^{k-s-1})^{-1}\cdot\sum_{n=1}^{\infty} a(n^2)\,n^{-s}.
$$
\end{remark}

\vskip 3mm
Let $A$ be a $g$-dimensional simple abelian variety defined over $\BQ$. For a prime $\ell$, we set
\begin{equation*}
  A[\ell^n]:=\{ \, x\in A(\overline{\BQ})\,|\ \ell^n\cdot x=0\,\}.
\end{equation*}
Then $A[\ell^n]\cong (\BZ/\ell^n \BZ)^g\times (\BZ/\ell^n \BZ)^g$\,(cf.\,\cite{MU}). Then the
Tate module of $A$ is given by
\begin{equation*}
  T_{\ell}(A):=\lim_{\longleftarrow} A[\ell^n]\cong \BZ_\ell^g\times \BZ_\ell^g
  \cong \BZ_{\ell}^{2g}.
\end{equation*}
Therefore we have the $2g$-dimensional $\ell$-adic Galois representation of $G_\BQ$
\begin{equation*}
\rho_{A,\ell}:G_\BQ\lrt GL(2g,\BZ_\ell)\subset GL(2g,\BQ_\ell).
\end{equation*}

\begin{definition}\label{def:3.4}
A $2g$-dimensional $\ell$-adic Galois representation $\rho$ of $G_\BQ$ given by
$$
\rho:G_\BQ\lrt GL(2g,\BZ_\ell)\subset GL(2g,\BQ_\ell)
$$
is called $\sf{modular}$ if there is a Siegel Hecke eigenform $F(\Om)\in [\G_{g,0}(N), g+1]_0$
of weight $g+1$ with respect to $\G_{g,0}(N)$ such that
\begin{equation*}
  {\rm Tr}\left( \rho({\rm Frob}_p)\right)=a(pI_g;F)\quad {\rm and}\quad
  \det \left( \rho({\rm Frob}_p)\right)=p^g\quad {\rm for\ any\ prime}\ p\not|\, \ell N,
\end{equation*}
where
\begin{equation*}
  F(\Om)=\sum_{T} a(T;F)\,e^{2\pi i\,{\rm Tr}(T\Om)}
\end{equation*}
is a Fourier expansion of $F(\Om)$.
\end{definition}

\begin{definition}\label{def:3.5}
Let $A$ be a $g$-dimensional simple abelian variety defined over $\BQ$ and let $\ell$ be a prime.
 For a prime $p$, we let
\begin{equation*}
  L_p(A,s):=\left\{ \det\left( I_{2g}-p^{-s}\cdot
  \rho_{A,\ell}({\rm Frob}_p)\big|_{T_{\ell}(A)}\right)\right\}^{-1}
\end{equation*}
be the local $L$-function of $A$ at $p$. We define the $L$-function $L(A,s)$ of $A$ by
\begin{equation*}
  L(A,s)=\prod_{p:\,{\rm prime}}L_p(A,s).
\end{equation*}
\end{definition}

\begin{definition}\label{def:3.6}
Let $A$ be a $g$-dimensional simple abelian variety defined over $\BQ$.
$A$ is called $\sf{modular}$ if there exists a Siegel Hecke eigenform
$F(\Om)\in [\G_{g,0}(N), g+1]_0$
of weight $g+1$ with respect to $\G_{g,0}(N)$ such that
\begin{equation*}
  L(A,s)=D(F,s),\ Z_F(s)\ {\rm or}\ D_F(s).
\end{equation*}
\end{definition}
\vskip 3mm
For two positive integers $g$ and $N$, we let
\begin{equation*}
 {\mathcal A}_{g,0}(N):=\G_{g,0}(N)\ba \BH_g
\end{equation*}
be the Siegel modular variety of level structure $N$ and let
${\mathcal A}_{g,0}^{\rm tor}(N)$ be a smooth toroidal compactification of
${\mathcal A}_{g,0}(N)$\,(cf.\,\cite{AMRT, F-C}). We denote by
\begin{equation*}
 \Omega^i \left( {\mathcal A}_{g,0}^{\rm tor}(N)\right), \quad 0\leq i\leq \frac{g(g+1)}{2}
\end{equation*}
the complex vector space of holomorphic $i$-forms on ${\mathcal A}_{g,0}^{\rm tor}(N)$.
The Jacobian variety ${\rm Jac}\left( {\mathcal A}_{g,0}^{\rm tor}(N) \right)$ of
${\mathcal A}_{g,0}^{\rm tor}(N)$
is defined to be the abelian variety
\begin{equation}\label{(3.2)}
 {{\rm Jac}}({\mathcal A}_{g,0}^{\rm tor}(N)):=\Omega^{\nu}
 \left( {\mathcal A}_{g,0}^{\rm tor}(N) \right)^{\vee}
 / H_\nu \left( {\mathcal A}_{g,0}^{\rm tor}(N),\BZ \right),\qquad \nu=\frac{g(g+1)}{2}.
\end{equation}
The geometric genus of ${\mathcal A}_{g,0}^{\rm tor}(N)$ is the dimension of
the Jacobian variety ${\rm Jac}\left( {\mathcal A}_{g,0}^{\rm tor}(N) \right)$.
It is known that the following two vector spaces are isomprphic\,:
\begin{equation}\label{(3.3)}
 [\G_{g,0}(N),g+1]_0\cong \Omega^{\nu}
 \left( {\mathcal A}_{g,0}^{\rm tor}(N) \right),\qquad \nu=\frac{g(g+1)}{2}.
\end{equation}
More precisely, for a coordinate $\Omega=(\omega_{ij})\in \BH_g$, we let
\begin{equation*}
  \omega_0:=d\omega_{11}\wedge d\omega_{12}\wedge d\omega_{13}\wedge \cdots \wedge
  d\omega_{gg}
\end{equation*}
be a holomorphic $\nu$-form on $\BH_g$. If $\omega=F(\Om)\,\omega_0$ is a $\G_{g,0}(N)$-invariant
holomorphic form on $\BH_g$,
then
\begin{equation*}
F(\g\cdot \Om)=\det (C\Om+D)^{g+1}F(\Om)
\end{equation*}
for all $\g=\begin{pmatrix}
              A & B \\
              C & D
            \end{pmatrix}\in \G_{g,0}(N)$ and $\Om\in\BH_g.$ Thus
$F\in [\G_{g,0}(N),g+1].$ It was shown by Freitag\,\cite{Fr} that $\om$ can be extended
to a holomorphic $\nu$-form on ${\mathcal A}_{g,0}^{\rm tor}(N)$ if and only if $F$ is
a cusp form in $[\G_{g,0}(N),g+1]_0.$ Indeed, the mapping
\begin{equation*}
[\G_{g,0}(N),g+1]_0 \lrt \Omega^{\nu}
 \left( {\mathcal A}_{g,0}^{\rm tor}(N) \right),\qquad F\longmapsto F\,\omega_0
\end{equation*}
is an isomorphism as complex vector spaces. We observe that if
$\om_k:=G(\Om)\, \om_0^{\otimes k}$ is a $\G_{g,0}(N)$-invariant
holomorphic form on $\BH_g$ of degree $k\nu$, then $G(\Om)\in [\G_{g,0}(N),k(g+1)]_0$
is a cusp form of weight $k(g+1).$
\vskip 3mm
Therefore according to \eqref{(3.2)} and \eqref{(3.3)}, we have
\begin{equation}\label{(3.4)}
{{\rm Jac}}({\mathcal A}_{g,0}^{\rm tor}(N))
 \cong [\G_{g,0}(N),g+1]_0^{\vee}
 / H_\nu \left( {\mathcal A}_{g,0}^{\rm tor}(N),\BZ \right),\qquad \nu=\frac{g(g+1)}{2}.
\end{equation}
If there is no confusion, we simply set
$$
J_{g,0}(N):= {{\rm Jac}}({\mathcal A}_{g,0}^{\rm tor}(N)).
$$

\vskip 3mm
We propose the following conjectures.
\vskip 2mm \noindent
\begin{conjecture}\label{conj:3.7}
A simple abelian variety of dimesion $g$ defined over $\BQ$ is modular.
\end{conjecture}
\vskip 3mm
\noindent
\begin{conjecture}\label{conj:3.8}
Let A be a simple abelian variety of dimesion $g$ defined over $\BQ$.
The following statements are equivalent\,:
\vskip 2mm
{\sf{(MAV1)}} $A$ is modular.
\vskip 1mm
{\sf{(MAV2)}} There exists a non-constant holomorphic mapping
${\mathcal A}_{g,0}^{\rm tor}(N)\lrt A$ for some
\par
\ \ \ \ \ \ \ \ \ \ \ positive integer $N$.
\vskip 1mm
{\sf{(MAV3)}} There exists a non-constant holomorphic mapping
$J_{g,0}(N)\lrt A$ for some
\par
\ \ \ \ \ \ \ \ \ \ \ positive integer $N$.
\vskip 1mm
{\sf{(MAV4)}} $\rho_{A,\ell}$ is modular for any prime $\ell$.
\end{conjecture}

\vskip 3mm
We propose the following problems.
\vskip 2mm \noindent
\begin{problem}\label{prob:3.9}
Let $F\in [\G_{g,0}(N),g+1]_0$ be a Siegel Hecke eigenform of weight $g+1$.
Associate to $F$ a $2g$-dimensional continuous irreducible Galois representation
of $G_\BQ$.
\end{problem}

\vskip 2mm \noindent
\begin{problem}\label{prob:3.10}
Let $k$ be a positive integer.
Let $F\in [\G_{g,0}(N),k]_0$ be a Siegel Hecke eigenform of weight $k$.
Associate to $F$ a $2g$-dimensional continuous irreducible Galois representation
of $G_\BQ$.
\end{problem}

\begin{remark}\label{rk:3.11}
As mentioned in Section 2, in the case $g=1$, to a Hecke eigenform of weight $k\geq 1$,
Shimura, Deligne and Serre \cite{Sh2, De, D-S} associated a two-dimensional continuous
irreducible Galois representation. For the case $g=2$, Taylor \cite{T1,T2} tried
to associate the four dimensional continuous Galois representation of $G_\BQ$ to
a Siegel Hecke eigenform of small weight. But he did not specify the precise weight.
\end{remark}

\end{section}

\vskip 10mm

\begin{section}{\bf Appendix : The Hecke Algebra and $p$-Satake Parameters}
\setcounter{equation}{0}
\vskip 0.5cm \noindent {\bf 4.1. The Structure of the Hecke Algebra}
\vskip 0.3cm For a positive
integer $g$, we let $\G_g=Sp(2g,\BZ)$ and let
\begin{equation*}
\Delta_g:=GSp(2g,\BQ)=\big\{\,M\in GL(2g,\BQ)\,|\
{}^tMJ_gM=l(M)J_g,\ l(M)\in \BQ^{\times}\,\big\}
\end{equation*}
be the group of symplectic similitudes of the rational symplectic
vector space $(\BQ^{2g},\langle\ ,\ \rangle)$. We put
\begin{equation*}
\Delta_g^+:=GSp(2g,\BQ)^+=\big\{\,M\in \Delta_g\,|\
 l(M)>0\,\big\}.
\end{equation*}

Following the notations in \cite{Fr}, we let ${\mathscr
H}(\G_g,\Delta_g)$ be the complex vector space of all formal
finite sums of double cosets $\G_g M\G_g$ with $M\in \Delta_g^+$.
A double coset $\G_g M\G_g\,(M\in \Delta_g^+)$ can be written as a
finite disjoint union of right cosets $\G_gM_{\nu}\,(1\leq \nu\leq
h)\,:$
\begin{equation*}
\G_g M\G_g= \bigcup^h_{\nu=1} \G_gM_{\nu}\quad ( \textrm{disjoint}).
\end{equation*}
Let ${\mathscr L}(\G_g,\Delta_g)$ be the complex vector space
consisting of formal finite sums of right cosets $\G_gM$ with
$M\in \Delta^+$. For each double coset $\G_g M\G_g= \bigcup^h_{\nu=1}
\G_gM_{\nu}$ we associate an element $j(\G_gM\G_g)$ in ${\mathscr
L}(\G_g,\Delta_g)$ defined by
\begin{equation*}
j(\G_g M\G_g):= \sum^h_{\nu=1} \G_gM_{\nu}.
\end{equation*}
Then $j$ induces a linear map
\begin{equation}\label{(4.1)}
j_*:{\mathscr H}(\G_g,\Delta_g)\lrt {\mathscr L}(\G_g,\Delta_g).
\end{equation}
We observe that $\Delta_g$ acts on ${\mathscr L}(\G_g,\Delta_g)$
as follows:
\begin{equation*}
\big( \sum_{j=1}^h c_j\,\G_gM_j\big)\cdot M=\sum_{j=1}^h c_j\,\G_g
M_jM,\quad M\in \Delta_g.
\end{equation*}
We denote
\begin{equation*}
{\mathscr L}(\G_g,\Delta_g)^{\G_g}:=\big\{\, T\in {\mathscr
L}(\G_g,\Delta_g)\,|\ T\cdot \g=T\ \textrm{for all}\ \g\in
\G_g\,\big\}
\end{equation*}
be the subspace of $\G_g$-invariants in ${\mathscr
L}(\G_g,\Delta_g)$. Then we can show that ${\mathscr
L}(\G_g,\Delta_g)^{\G_g}$ coincides with the image of $j_*$ and
the map
\begin{equation}\label{(4.2)}
j_*:{\mathscr H}(\G_g,\Delta_g)\lrt {\mathscr
L}(\G_g,\Delta_g)^{\G_g}
\end{equation}
is an isomorphism of complex vector spaces
(cf.\,\cite[p.\,228]{Fr}). From now on we identify ${\mathscr
H}(\G_g,\Delta_g)$ with ${\mathscr L}(\G_g,\Delta_g)^{\G_g}.$

We define the multiplication of the double coset $\G_g M\G_g$ and
$\G_g N$ by
\begin{equation}\label{(4.3)}
(\G_g M\G_g)\cdot(\G_g N)=\sum_{j=1}^h \G_gM_jN,\quad M,N\in
\Delta_g,
\end{equation}
where $\G_g M\G_g= \bigcup^h_{j=1} \G_gM_{j}\ ( \textrm{disjoint}).$
The definition (4.3) is well defined, i.e., independent of the
choice of $M_j$ and $N$. We extend this multiplication to
${\mathscr H}(\G_g,\Delta_g)$ and ${\mathscr L}(\G_g,\Delta_g)$.
Since
\begin{equation*}
{\mathscr H}(\G_g,\Delta_g)\cdot {\mathscr
H}(\G_g,\Delta_g)\subset {\mathscr H}(\G_g,\Delta_g),
\end{equation*}
${\mathscr H}(\G_g,\Delta_g)$ is an associative algebra with the
identity element $\G_gI_{2g}\G_g=\G_g$. The algebra ${\mathscr
H}(\G_g,\Delta_g)$ is called the $ \textit{Hecke algebra}$ with
respect to $\G_g$ and $\Delta_g$. \vskip 0.2cm We now describe the
structure of the Hecke algebra ${\mathscr H}(\G_g,\Delta_g)$. For
a prime $p$, we let $\BZ[1/p]$ be the ring of all rational numbers
of the form $a\cdot p^{\nu}$ with $a,\nu\in\BZ.$ For a prime $p$,
we denote
\begin{equation*}
\Delta_{g,p}:=\Delta_g \cap GL\big(2g,\BZ[1/p]\big).
\end{equation*}
Then we have a decomposition of ${\mathscr H}(\G_g,\Delta_g)$
\begin{equation*}
{\mathscr H}(\G_g,\Delta_g)=\bigotimes_{p\,:\, \textrm{prime}}
{\mathscr H}(\G_g,\Delta_{g,p})
\end{equation*}
as a tensor product of local Hecke algebras ${\mathscr
H}(\G_g,\Delta_{g,p}).$ We denote by $\check{\mathscr
H}(\G_g,\Delta_g)$\,(resp. $\check{\mathscr H}(\G_g,\Delta_{g,p})$
the subring of ${\mathscr H}(\G_g,\Delta_g)$\,(resp. ${\mathscr
H}(\G_g,\Delta_{g,p})$ by integral matrices.

In order to describe the structure of local Hecke operators
${\mathscr H}(\G_g,\Delta_{g,p})$, we need the following lemmas.

\begin{lemma}\label{lem:4.1}
Let $M\in \Delta_g^+$ with ${}^tMJ_gM=lJ_g$. Then the
double coset $\G_gM\G_g$ has a unique representative of the form
\begin{equation*}
M_0= {\rm diag}(a_1,\cdots,a_g,d_1,\cdots,d_g),
\end{equation*}
\noindent where $a_g|d_g,\ a_j>0,\ a_jd_j=l$ for $1\leq j\leq g$
and $a_k| a_{k+1}$ for $1\leq k\leq g-1.$
\end{lemma}

For a positive integer $l$, we let
\begin{equation*}
O_g(l):=\big\{\, M\in GL(2g,\BZ)\,|\ {}^tMJ_gM=lJ_g\ \big\}.
\end{equation*}
Then we see that $O_g(l)$ can be written as a finite disjoint
union of double cosets and hence as a finite union of right
cosets. We define $T(l)$ as the element of ${\mathscr H}(\G_g,\Delta_g)$
defined by $O_g(l).$

\begin{lemma}\label{lem:4.2}
(a) Let $l$ be a positive integer. Let
\begin{equation*}
O_g(l)=\cup_{\nu=1}^h\G_gM_{\nu}\quad ( {\rm disjoint})
\end{equation*}
be a disjoint union of right cosets $\G_gM_{\nu}\,(1\leq \nu\leq
h).$ Then each right coset $\G_gM_{\nu}$ has a representative of
the form
\begin{equation*}
M_{\nu}=\begin{pmatrix} A_{\nu}&B_{\nu}\\
                   0&   D_{\nu}\end{pmatrix}, \quad
                   {}^tA_{\nu}D_{\nu}=l\,I_g,\quad A_{\nu}\
                   {\rm is\ upper\ triangular}.
\end{equation*}
(b) Let $p$ be a prime. Then
\begin{equation*}
T(p)=O_g(p)=\G_g\begin{pmatrix} I_g&0\\
                   0&  pI_g\end{pmatrix}\G_g
\end{equation*}
and
\begin{equation*}
T(p^2)=\sum_{i=0}^g T_i(p^2),
\end{equation*}
where
\begin{equation*}
T_k(p^2):=\begin{pmatrix} I_{g-k}&0 & 0 & 0\\
                   0&  pI_k  & 0 & 0 \\ 0 & 0 & p^2I_{g-k} & 0\\
                   0 & 0 & 0 & pI_k\end{pmatrix}\G_g,\quad 0\leq
                   k\leq g.
\end{equation*}
\end{lemma}
\noindent $ \textit{Proof.}$ The proof can be found in
\cite[p.\,225 and p.\,250]{Fr}.
\hfill$\square$

\vskip 3mm
For example, $T_g(p^2)=\G_g(pI_{2g})\G_g$ and
\begin{equation*}
T_0(p^2)=\G_g\begin{pmatrix} I_g&0\\
                   0&  p^2I_g\end{pmatrix}\G_g=T(p)^2.
\end{equation*}

We have the following
\begin{theorem}\label{thm:4.3}
The local Hecke algebra $\check{\mathscr H}(\G_g,\Delta_{g,p})$
is generated by algebraically independent generators
$T(p),\,T_1(p^2),\cdots,T_g(p^2).$
\end{theorem}
\noindent $ \textit{Proof.}$ The proof can be found in
\cite[p.\,250 and p.\,261]{Fr}. \hfill$\square$

\vskip 0.2cm On $\D_g$ we have the anti-automorphism $M\mapsto
M^*:=l(M)M^{-1}\,(M\in \D_g)$. Obviously $\G_g^*=\G_g$. By Lemma
4.1, $(\G_gM\G_g)^*=\G_gM^*\G_g=\G_gM\G_g.$ According to
\cite{Sh2}, Proposition 3.8, ${\mathscr{H}}(\G_g,\Delta_{g})$ is
commutative.
\vskip 0.2cm Let $X_0,X_1,\cdots,X_g$ be the $g+1$
variables. We define the automorphisms
\begin{equation*}
w_j: \BC\big[X_0^{\pm 1},X_1^{\pm 1},\cdots,X_g^{\pm 1}\big]\lrt
\BC\big[X_0^{\pm 1},X_1^{\pm 1},\cdots,X_g^{\pm 1}\big],\quad
1\leq j\leq g
\end{equation*}
by
\begin{equation*}
w_j(X_0)=X_0X_j^{-1},\ \ \ w_j(X_j)=X_j^{-1},\ \ \ w_j(X_k)=X_k\ \
\textrm{for} \ k\neq 0,j.
\end{equation*}
Let $W_g$ be the finite group generated by $w_1,\cdots,w_g$ and
the permutations of variables $X_1,\cdots,X_g$. Obviously $w_j^2$
is the identity map and $|W_g|=2^g g!$.

\begin{theorem}\label{thm:4.4}
There exists an isomorphism
\begin{equation*}
Q:{\mathscr H}(\G_g,\Delta_{g,p})\lrt \BC\big[X_0^{\pm 1
},X_1^{\pm 1 },\cdots,X_g^{\pm 1}\big]^{W_g}.
\end{equation*}
In fact, $Q$ is defined by
\begin{equation*}
Q\big( \sum_{j=1}^h\G_gM_j\big)=\sum_{j=1}^h
Q(\G_gM_j)=\sum_{j=1}^h
X_0^{-k_0(j)}\prod_{\nu=1}^g\big(p^{-\nu}X^{\nu}\big)^{k_{\nu}(j)}|\det
A_j|^{g+1},
\end{equation*}
where we chose the representative $M_j$ of $\G_g M_j$ of the form
\begin{equation*}
M_j=\begin{pmatrix} A_j& B_j\\
                   0&  D_j\end{pmatrix},\quad A_j=\begin{pmatrix} p^{k_1(j)}& \ldots & *\\
                   0 & \ddots & \vdots \\
                   0&  0 & p^{k_g(j)}\end{pmatrix}.
\end{equation*}
We note that the integers $k_1(j),\cdots,k_g(j)$ are uniquely
determined.
\end{theorem}
\noindent $ \textit{Proof.}$ The proof can be found in \cite[pp.\,254--261]{Fr}.
\hfill $\square$

\vskip 0.3cm For a prime $p$, we let
\begin{equation*}
{\mathscr H}(\G_g,\Delta_{g,p})_{\BQ}:=\left\{\, \sum c_j\,\G_g
M_j\G_g \in {\mathscr H}(\G_g,\Delta_{g,p})\,|\ c_j\in\BQ\
\right\}
\end{equation*}
be the $\BQ$-algebra contained in ${\mathscr
H}(\G_g,\Delta_{g,p})$. We put
\begin{equation*}
G_p:=GSp(g,\BQ_p)\quad \textrm{and}\quad K_p=GSp(g,\BZ_p).
\end{equation*}
We can identify ${\mathscr H}(\G_g,\Delta_{g,p})_{\BQ}$ with the
$\BQ$-algebra ${\mathscr H}_{g,p}^{\BQ}$ of $\BQ$-valued locally
constant, $K_p$-biinvariant functions on $G_p$ with compact
support. The multiplication on ${\mathscr H}_{g,p}^{\BQ}$ is given
by
\begin{equation*}
(f_1*f_2)(h)=\int_{G_p}f_1(g)\,f_2(g^{-1}h)dg,\quad f_1,f_2\in
{\mathscr H}_{g,p}^{\BQ},
\end{equation*}
where $dg$ is the unique Haar measure on $G_p$ such that the
volume of $K$ is $1$. The correspondence is obtained by sending
the double coset $\G_gM\G_g$ to the characteristic function of
$K_pMK_p$. \vskip 0.2cm In order to describe the structure of
${\mathscr H}_{g,p}^{\BQ}$, we need to understand the $p$-adic
Hecke algebras of the diagonal torus ${\mathbb T}$ and the Levi
subgroup ${\mathbb M}$ of the standard parabolic group. Indeed,
${\mathbb T}$ is defined to be the subgroup consisting of diagonal
matrices in $\D_g$ and
\begin{equation*}
{\mathbb M}=\left\{\,\begin{pmatrix} A& 0\\
                   0&  D\end{pmatrix}\in \D_g\ \right\}
\end{equation*}
is the Levi subgroup parabolic subgroup
\begin{equation*}
\left\{\,\begin{pmatrix} A& B\\
                   0&  D\end{pmatrix}\in \D_g\ \right\}.
\end{equation*}
Let $Y$ be the co-character group of ${\mathbb T}$, i.e., $Y=
\textrm{Hom}({\mathbb G}_m,{\mathbb T}).$ We define the local
Hecke algebra ${\mathscr H}_p({\mathbb T})$ for ${\mathbb T}$ to
be the $\BQ$-algebra of $\BQ$-valued, ${\mathbb
T}(\BZ_p)$-biinvariant functions on $\BT(\BQ_p)$ with compact
support. Then ${\mathscr H}_p(\BT)\cong \BQ[Y],$ where $\BQ[Y]$ is
the group algebra over $\BQ$ of $Y$. An element $\la\in Y$
corresponds the characteristic function of the double coset
$D_{\la}=K_p\la(p)K_p$. It is known that  ${\mathscr H}_p({\mathbb
T})$ is isomorphic to the ring $\BQ\big[(u_1/v_1)^{\pm
1},\cdots,(u_g/v_g)^{\pm 1},(v_1\cdots v_g)^{\pm 1}\big]$ under
the map
\begin{equation*}
(a_1,\cdots,a_g,c)\mapsto (u_1/v_1)^{a_1}\cdots
(u_g/v_g)^{a_g}(v_1\cdots v_g)^{c}.
\end{equation*}
Similarly we have a $p$-adic Hecke algebra ${\mathscr
H}_p({\mathbb M})$. Let $W_{\D_g}=N(\BT)/\BT$ be the Weyl group
with respect to $(\BT,\D_g)$, where $N(\BT)$ is the normalizer of
$\BT$ in $\D_g$. Then $W_{\D_g}\cong S_g \ltimes (\BZ/2\BZ)^g,$
where the generator of the $i$-th factor $\BZ/2\BZ$ acts on a
matrix of the form $ \textrm{diag}(a_1,\cdots,a_g,d_1,\cdots,d_g)$
by interchanging $a_i$ and $d_i$, and the symmetry group $S_g$
acts by permuting the $a_i${'}$s$ and $d_i${'}$s$. We note that
$W_{\D_g}$ is isomorphic to $W_g$. The Weyl group $W_{\BM}$ with
respect to $(\BT,\BM)$ is isomorphic to $S_g$. We can prove that
the algebra ${\mathscr H}_p({\mathbb T})^{W_{\D_g}}$ of
$W_{\D_g}$-invariants in ${\mathscr H}_p({\mathbb T})$ is
isomorphic to $\BQ \big[ Y_0^{\pm 1},Y_1,\cdots,
Y_g\big]$\,(cf.\,\cite{Fr}). We let
\begin{equation*}
B=\left\{\,\begin{pmatrix} A& B\\
                   0&  D\end{pmatrix}\in \D_g\,\Big|\ A\
                   \textrm{is upper triangular,}\ D\ \textrm{is lower trianular}\ \right\}
\end{equation*}
be the Borel subgroup of $\D_g$. A set $\Phi^+$ of positive roots
in the root system $\Phi$ determined by $B$. We set $\rho={\frac
12}\sum_{\alpha \in \Phi^+}\alpha.$

\vskip 0.2cm Now we have the map $\alpha_{\BM}:\BM\lrt {\mathbb
G}_m$ defined by
\begin{equation*}
\alpha_{\BM}(M):=l(M)^{-{{g(g+1)}\over 2}} \big( \det
A\big)^{g+1},\quad M=\begin{pmatrix}A & 0\\ 0 & D\end{pmatrix}\in
\BM
\end{equation*}
and the map $\beta_{\BT}:\BT\lrt {\mathbb G}_m$ defined by
\begin{equation*}
\beta_{\BT}(\textrm{diag}(a_1,\cdots,a_g,d_1,\cdots,d_g)):=\prod_{i=1}^g
a_1^{g+1-2i},\quad \textrm{diag}(a_1,\cdots,a_g,d_1,\cdots,d_g)\in
\BT.
\end{equation*}
Let $\theta_{\BT}:=\alpha_{\BM}\,\beta_{\BT}$ be the character of
$\BT.$ The $ \textit{Satake's spherical map}\ S_{p,\BM}:{\mathscr
H}_{g,p}^{\BQ}\lrt {\mathscr H}_p(\BM)$ is defined by
\begin{equation}\label{(4.4)}
S_{p,\BM}(\phi)(m):=|\alpha_\BM(m)|_p
\int_{U(\BQ_p)}\phi(mu)du,\quad \phi\in {\mathscr H}_{g,p}^{\BQ},\
m\in \BM,
\end{equation}
where $|\ \ |_p$ is the $p$-adic norm and $U(\BQ_p)$ denotes the
unipotent radical of $\D_g$. Also another $\textit{Satake's
spherical map}\ S_{\BM,\BT}:{\mathscr H}_p({\BM})\lrt {\mathscr
H}_p(\BT)$ is defined by
\begin{equation}\label{(4.5)}
S_{\BM,\BT}(f)(t):=|\beta_\BT(t)|_p \int_{\BM\cap{\mathbb
N}}f(tn)dn,\quad t\in {\mathscr H}_p(\BT),\ t\in \BT,
\end{equation}
where ${\mathbb N}$ is a nilpotent subgroup of $\D_g$.

\begin{theorem}\label{thm:4.5}
The Satake's spherical maps $S_{p,\BM}$ and
$S_{\BM,\BT}$ define the isomorphisms of $\BQ$-algebras
\begin{equation}\label{(4.6)}
{\mathscr H}_{g,p}^{\BQ}\cong {\mathscr H}_p(\BT)^{W_{\D_g}} \quad
\textrm{and}\quad {\mathscr H}_p(\BM)\cong {\mathscr
H}_p(\BT)^{W_{\BM}}.
\end{equation}
\end{theorem}
\noindent We define the elements $\phi_k\,(0\leq k\leq g)$ in
${\mathscr H}_p(\BM)$ by
\begin{equation*}
\phi_k:=p^{-{{k(k+1)}\over 2}}\,\BM(\BZ_p)\begin{pmatrix} I_{g-k}
& 0 & 0 \\ 0 & pI_g & 0 \\ 0 & 0 & I_k \end{pmatrix}
\BM(\BZ_p),\quad i=0,1,\cdots,g.
\end{equation*}
Then we have the relation
\begin{equation}\label{(4.7)}
S_{p,\BM}(T(p))=\sum_{k=0}^g\phi_k
\end{equation}
and
\begin{equation}\label{(4.8)}
S_{p,\BM}\big(T_i(p^2)\big)=\sum_{j,k\geq 0,\,i+j\leq k}
m_{k-j}(i)\, p^{- {k-j+1\choose 2}}\phi_j\phi_k,
\end{equation}
where
\begin{equation*}
m_s(i):=\sharp \left\{\,A\in M(s,{\mathbb F}_p)\,|\ {}^tA=A,\quad
\textrm{corank}(A)=i\ \right\}.
\end{equation*}
Moreover, for $k=0,1,\cdots,g$, we have
\begin{equation}\label{(4.9)}
S_{\BM,\BT}(\phi_k)=(v_1\cdots v_g)E_k(u_1/v_1,\cdots,u_g/v_g),
\end{equation}
where $E_k$ denotes the elementary symmetric function of degree
$k$. The proof of (4.7)--(4.9) can be found in
\cite[pp.\,142--145]{A2}.

\vskip 0.5cm \noindent
{\bf 4.2. Action of the Hecke Algebra on Siegel Modular Forms}

\vskip 0.3cm
Let $(\rho,V_{\rho})$ be a finite dimensional
irreducible representation of $GL(g,\BC)$ with highest weight
$(k_1,\cdots,k_g)$. For a function $F:\BH_g\lrt V_{\rho}$ and
$M\in \D_g^+,$ we define
\begin{equation*}
(f|_{\rho}M)(\Om)=\rho(C\Omega+D)^{-1}f(M\cdot\Om),\quad
M=\begin{pmatrix}A & B\\ C & D \end{pmatrix}\in \D_g^+.
\end{equation*}
It is easily checked that $f|_{\rho}M_1M_2=\big(
f|_{\rho}M_1\big)|_{\rho}M_2$ for $M_1,M_2\in \D_g^+.$

\vskip 0.2cm
We now consider a subset ${\mathscr M}$ of $\D_g$
satisfying the following properties (M1) and (M2)\,:
\par \ \ (M1)\ \ \ ${\mathscr M}=\bigcup_{j=1}^h \G_gM_j\quad$ (disjoint
union);
\par \ \ (M2)\ \ \ ${\mathscr M}\,\G_g\subset {\mathscr M}.$

\vskip 0.2cm
For a Siegel modular form $f\in M_{\rho}(\G_g)$, we define
\begin{equation}\label{(4.10)}
T({\mathscr M})f:=\sum_{j=1}^h f|_{\rho}M_j.
\end{equation}
This is well defined, i.e., is independent of the choice of
representatives $M_j$ because of the condition (M1). On the other
hand, it follows from the condition (M2) that $T({\mathscr
M})f|_{\rho}\g=T({\mathscr M})f$ for all $\g\in \G_g.$ Thus we get
a linear operator
\begin{equation}\label{(4.11)}
T({\mathscr M}):M_{\rho}(\G_g)\lrt M_{\rho}(\G_g).
\end{equation}
We know that each double coset $\G_gM\G_g$ with $M\in \D_g$
satisfies the condition $(M1)$ and $(M2)$. Thus a linear operator
$T({\mathscr M})$ defined in Formula (4.10) induces naturally the action of
the Hecke algebra ${\mathscr H}(\G_g,\Delta_g)$ on
$M_{\rho}(\G_g)$. More precisely, if ${\mathscr N}=\sum_{j=1}^h
c_j\G_gM_j\G_g\in {\mathscr H}(\G_g,\Delta_g)$, we define
\begin{equation*}
T({\mathscr N})=\sum_{j=1}^h c_jT(\G_gM_j\G_g).
\end{equation*}
Then $T({\mathscr N})$ is an endomorphism of $M_{\rho}(\G_g)$.

\vskip 0.2cm
Now we fix a Siegel modular form $F$ in
$M_{\rho}(\G_g)$ which is an eigenform of the Hecke algebra
${\mathscr H}(\G_g,\Delta_g)$. Then we obtain an algebra
homomorphism $\la_F:{\mathscr H}(\G_g,\Delta_g)\lrt \BC$
determined by
\begin{equation*}
T(F)=\la_F(T)F,\quad T\in {\mathscr H}(\G_g,\Delta_g).
\end{equation*}
By Theorem 4.4 or Theorem 4.5, one has
\begin{eqnarray*}
{\mathscr H}(\G_g,\Delta_{g,p})&\cong& {\mathscr
H}_{g,p}^{\BQ}\otimes \BC \cong \BC[Y]^{W_g}\\
&\cong& {\mathscr H}_p({\mathbb T})^{W_g}\otimes \BC\\
&\cong& \BC\big[(u_1/v_1)^{\pm 1},\cdots,(u_g/v_g)^{\pm
1},(v_1\cdots v_g)^{\pm 1}\big]^{W_g}\\
&\cong& \BC[Y_0,Y_0^{-1},Y_1,\cdots,Y_g],
\end{eqnarray*}
where $Y_0,Y_1,\cdots,Y_g$ are algebraically independent.
Therefore one obtains an isomorphism
\begin{equation*}
\textrm{Hom}_{\BC}\big({\mathscr
H}(\G_g,\Delta_{g,p}),\BC\big)\cong
\textrm{Hom}_{\BC}\big({\mathscr H}_{g,p}^{\BQ}\otimes
\BC,\BC\big)\cong (\BC^{\times})^{(g+1)}/W_g.
\end{equation*}
The algebra homomorphism $\la_F\in
\textrm{Hom}_{\BC}\big({\mathscr H}(\G_g,\Delta_{g,p}),\BC\big)$
is determined by the $W_g$-orbit of a certain $(g+1)$-tuple $\big(
\alpha_{F,0},\alpha_{F,1},\cdots,\alpha_{F,g}\big)$ of nonzero
complex numbers, called the $p$-$ \textit{Satake parameters}$ of
$F$. For brevity, we put $\alpha_i=\alpha_{F,i},\ i=0,1,\cdots,g$.
Therefore $\alpha_i$ is the image of $u_i/v_i$ and $\alpha_0$ is
the image of $v_1\cdots v_g$ under the map $\Theta$. Each
generator $w_i\in W_{\D_g}\cong W_g$ acts by
\begin{equation*}
w_j(\alpha_0)=\alpha_0\alpha_j^{-1}\quad
w_j(\alpha_j)=\alpha_j^{-1},\quad w_j(\alpha_k)=0\ \textrm{if}\
k\neq 0,j.
\end{equation*}
These $p$-Satake parameters $\alpha_0,\alpha_1\cdots,\alpha_g$
satisfy the relation
\begin{equation*}
\alpha_0^2\alpha_1\cdots \alpha_g=p^{\sum_{i=1}^gk_i-g(g+1)/2}.
\end{equation*}
Formula (4.12) follows from the fact that
$T_g(p^2)=\G_g(pI_{2g})\G_g$ is mapped to
$$p^{-g(g+1)/2}\,(v_1\cdots
v_g)^2\prod_{i=1}^g (u_i/v_i).$$ We refer to \cite[p.\,258]{Fr}
for more detail. According to Formulas (4.7)--(4.9), the eigenvalues
$\la_F\big(T(p)\big)$ and $\la_F\big(T_i(p^2)\big)$ with $1\leq
i\leq g$ are given respectively by
\begin{equation}\label{(4.12)}
\la_F\big(T(p)\big)=\alpha_0(1+E_1+E_2+\cdots+E_g)
\end{equation}
and
\begin{equation}\label{(4.13)}
\la_F\big(T_i(p^2)\big)=\sum_{j,k\geq 0,\,j+i\leq k}^g
m_{k-j}(i)\,p^{-{ {k-j+1}\choose 2}}\,\alpha_0^2E_jE_k,\quad
i=1,\cdots,g,
\end{equation}
where $E_j$ denotes the elementary symmetric function of degree
$j$ in the variables $\alpha_1,\cdots,\alpha_g$. The point is that
the above eigenvalues $\la_F\big(T(p)\big)$ and
$\la_F\big(T_i(p^2)\big)\ (1\leq i\leq g)$ are described in terms
of the $p$-Satake parameters $\alpha_0,\alpha_1\cdots,\alpha_g$.

\vskip 0.2cm \noindent
\begin{example}
Suppose
$g(\tau)=\sum_{n\geq 1}a(n)\,e^{2\pi i n\tau}$ is a normalized
eigenform in $S_k(\G_1)$. Let $p$ be a prime. Let $\beta$ be a
complex number determined by the relation
\begin{equation*}
(1-\beta X)(1-{\bar \beta}X)=1-a(p)X+p^{k-1}X^2.
\end{equation*}
Then
\begin{equation*}
\beta+{\bar\beta}=a(p)\quad \textrm{and}\quad \beta {\bar
\beta}=p^{k-1}.
\end{equation*}
The $p$-Satake parameters $\alpha_0$ and $\alpha_1$ are given by
\begin{equation*}
(\alpha_0,\alpha_1)=\left( \beta, {{\bar\beta}\over
\beta}\right)\quad or \quad \left( {\bar\beta}, {{\beta}\over
{\bar\beta}}\right).
\end{equation*}
It is easily checked that
$\alpha_0^2\alpha_1=\beta{\bar\beta}=p^{k-1}$ (cf. Formula
(4.12)).
\end{example}
\vskip 0.2cm\noindent
\begin{example}
For a positive integer $k$ with
$k>g+1$, we let
\begin{equation*}
G_k(\Om):=\sum_{M\in \G_{g,0}\backslash \G_g} \det
(C\Om+D)^k,\quad M=\begin{pmatrix} A & B\\ C & D\end{pmatrix}
\end{equation*}
be the Siegel Eisenstein series of weight $k$ in $M_k(\G_g)$,
where
$$\G_{g,0}:=\left\{ \begin{pmatrix} A & B\\ 0 &
D\end{pmatrix}\in \G_g\right\}$$ is a parabolic subgroup of
$\G_g$. It is known that $G_k$ is an eigenform of all the Hecke
operators\,(cf.\,\cite[p.\,268]{Fr}). Let $S_1,\cdots,S_h$ be a
complete system of representatives of positive definite even
unimodular integral matrices of degree $2k$.
We define the theta series $\theta_{S_{\nu}}$ by
\begin{equation*}
  \theta_{S_\nu}(\Omega):=\sum_{A\in \BZ^{(2k,g)}}
  e^{\pi\,i\,\rm{Tr}(S_\nu [A]\Omega)},\quad 1\leq \nu\leq h.
\end{equation*}
If $k>g+1$, the
Eisenstein series $G_k$ can be expressed as the weighted mean of
theta series $\theta_{S_1},\cdots,\theta_{S_h}\,$:
\begin{equation}
G_k(\Om)=\sum_{\nu=1}^h m_{\nu}\,\theta_{S_{\nu}}(\Om),\quad
\Om\in \BH_g,
\end{equation}
where
\begin{equation*}
m_{\nu}={ {A(S_{\nu},S_{\nu})^{-1}}\over
{A(S_1,S_1)^{-1}+\cdots+A(S_h,S_h)^{-1}} },\quad 1\leq \nu\leq h.
\end{equation*}
We recall that for two symmetric integral matrices $S$ of
degree $m$ and $T$ of degree $n$, $A(S,T)$ is defined by
\begin{equation*}
A(S,T):=\sharp \big\{\, G\in \BZ^{(m,n)}\,|\
S[G]=\,{}^tGSG=T\,\big\}.
\end{equation*}
Formula (4.14) was obtained by Witt \cite{Wit} as a special case
of the analytic version of Siegel's Hauptsatz.
\end{example}

\end{section}

\end{document}